\documentclass[11pt,reqno]{amsart}

\usepackage{amsthm}
\usepackage{amscd}
\usepackage{amsfonts}
\usepackage{amssymb}
\usepackage{amsgen}
\usepackage{amsmath}
\usepackage{amsopn}
\usepackage{verbatim}
\usepackage{xypic}
\usepackage{xspace}
\usepackage{multicol}
\usepackage{url}
\usepackage{upref}
\usepackage{pgf}
\usepackage{tikz}
\usetikzlibrary{arrows}

\theoremstyle{plain}
\newtheorem{thm}{Theorem}[section]
\newtheorem{lem}{Lemma}[section]
\newtheorem{prop}{Proposition}[section]

\theoremstyle{definition}
\newtheorem*{defn}{Definition}
\newtheorem*{rem}{Remark}

\numberwithin{equation}{section}

 \font\cyr=wncyr10
 \newcommand{\nc}{\newcommand}

\nc{\per}[1]{\underset{#1}{\boldsymbol \pi}\,}

 \nc{\MT}{{\rm MT}}
  \nc{\bgz}{{\bar{\gz}}}
 \nc{\wt}{{\rm wt}}
 \nc{\wht}{{\widehat}}
 \nc{\bwg}{{\bigwedge}}
 \nc{\mmu}{{\boldsymbol{\mu}}}
 \nc{\mal}{{{\scriptstyle \maltese}}}
 \nc{\fA}{{\mathfrak A}}
 \nc{\HH}{{\mathfrak H}}
 \nc{\ra}{\rightarrow}
 \nc{\ors}{{\vec s\,}}
 \nc{\os}{{\overset}}
 \nc{\G}{{\mathbb G}}
 \nc{\Z}{{\mathbb Z}}
 \nc{\R}{{\mathbb R}}
 \nc{\N}{{\mathbb N}}
 \nc{\ZN}{{\mathbb Z_{\ge 0}}}
 \nc{\Q}{{\mathbb Q}}
 \nc{\C}{{\mathbb C}}
 \nc{\Cnn}{{\mathbb C}_{\ge 0}}
 \nc{\Cp}{{\mathbb C}_{>0}}
 \nc{\MPV}{{\mathcal{MPV}}}
 
 \nc{\tB}{{\tilde B}}
 \nc{\tI}{{\tilde I}}
 \nc{\tJ}{{\tilde J}}
 \nc{\tK}{{\tilde K}}
 \nc{\Li}{{\rm Li}}
 \nc{\suf}{{\ast\,}}
 \nc{\sufq}{{\ast_q\,}}
 \nc{\gam}{{\gamma}}
 \nc{\gG}{{\Gamma}}
 \nc{\om}{{\omega}}
 \nc{\vep}{{\varepsilon}}
 \nc{\ga}{{\alpha}}
 \nc{\gl}{{\lambda}}
 \nc{\gb}{{\beta}}
 \nc{\gd}{{\delta}}
 \nc{\gs}{{\sigma}}
 \nc{\gS}{{\Sigma}}

 \nc{\gk}{{\kappa}}
 \nc{\tgz}{{\tilde{\zeta}}}
 \nc{\gO}{{\Omega}}
 \nc{\sif}{{\mathcal S}}
 \nc{\gt}{{\tau}}
\nc{\gz}{{\zeta}}
 \nc{\Lra}{\Longrightarrow}
 \nc{\lra}{\longrightarrow}
 \nc{\fS}{{\mathfrak S}}
 \nc{\DD}{{\mathfrak D}}

 \nc{\Llra}{\Longleftrightarrow}
 \nc{\ol}{\overline}
 \nc{\lms}{\longmapsto}
 \nc{\cv}{{{\mathsf c}{\mathsf v}}}
 \nc{\zq}{{\zeta_q}}
 \nc\qup{{q\uparrow 1}}
 \nc{\us}{\underset}
 \nc{\tn}{{\tilde{n}}}
 \nc{\gD}{{\Delta}}
 \nc{\bi}{{\bf i}}
 \nc{\bfone}{{\bf 1}}
 \nc{\bfa}{{\bf a}}
 \nc{\bfb}{{\bf b}}
 \nc{\bfc}{{\bf c}}
 \nc{\bfd}{{\bf d}}
 \nc{\bfe}{{\bf e}}
 \nc{\bff}{{\bf f}}
 \nc{\bfg}{{\bf g}}
 \nc{\bfh}{{\bf h}}
 \nc{\bfi}{{\bf i}}
 \nc{\bfj}{{\bf j}}
 \nc{\bfn}{{\bf n}}
 \nc{\bfl}{{\bf l}}
 \nc{\bfk}{{\bf k}}
 \nc{\bfm}{{\bf m}}
 \nc{\bfo}{{\bf o}}
 \nc{\bfp}{{\bf p}}
 \nc{\bfq}{{\bf q}}
 \nc{\bfr}{{\bf r}}
 \nc{\tbfs}{{\tilde{\bf s}}}
 \nc{\bfs}{{\bf s}}
 \nc{\hbfs}{{\hat{\bf s}}}
 \nc{\hs}{{\hat{s}}}
 \nc{\ts}{\tilde{s}}
 \nc{\bft}{{\bf t}}
 \nc{\bfu}{{\bf u}}
 \nc{\bfv}{{\bf v}}
 \nc{\bfw}{{\bf w}}
 \nc{\bfx}{{\bf x}}
 \nc{\bfy}{{\bf y}}
 \nc{\bfz}{{\bf z}}
 \nc{\bfB}{{\bf B}}
 \nc{\bfP}{{\bf P}}
 \nc{\bfQ}{{\bf Q}}
 \nc{\bfY}{{\bf Y}}
 \nc{\bfgb}{{\boldsymbol \gb}}
 \nc{\bfga}{{\boldsymbol \ga}}
 \nc{\bfrho}{{\boldsymbol \rho}}
 \nc{\bfchi}{{\boldsymbol \chi}}
 \nc{\QX}{{\Q\langle \bfX\rangle}}
 \nc{\QY}{{\Q\langle \bfY\rangle}}
 \nc{\CX}{{\C\langle \bfX\rangle}}
 \nc{\CY}{{\C\langle \bfY\rangle}}
 \nc{\QXX}{{\Q\langle\!\langle \bfX\rangle\!\rangle}}
 \nc{\QYY}{{\Q\langle\!\langle \bfY\rangle\!\rangle}}
 \nc{\CXX}{{\C\langle\!\langle \bfX\rangle\!\rangle}}
 \nc{\CYY}{{\C\langle\!\langle \bfY\rangle\!\rangle}}

 \nc{\bbA}{{\mathbb A}}
 \nc{\bbB}{{\mathbb B}}
 \nc{\bbC}{{\mathbb C}}
 \nc{\bbD}{{\mathbb D}}
 \nc{\bbE}{{\mathbb E}}
 \nc{\bbF}{{\mathbb F}}
 \nc{\bbG}{{\mathbb G}}
 \nc{\bbH}{{\mathbb H}}
 \nc{\bbI}{{\mathbb I}}
 \nc{\bbJ}{{\mathbb J}}
 \nc{\bbK}{{\mathbb K}}
 \nc{\bbL}{{\mathbb L}}
 \nc{\bbM}{{\mathbb M}}
 \nc{\bbN}{{\mathbb N}}
 \nc{\bbO}{{\mathbb O}}
 \nc{\bbP}{{\mathbb P}}
 \nc{\bbQ}{{\mathbb Q}}
 \nc{\bbR}{{\mathbb R}}
 \nc{\bbS}{{\mathbb S}}
 \nc{\bbT}{{\mathbb T}}
 \nc{\bbU}{{\mathbb U}}
 \nc{\bbV}{{\mathbb V}}
 \nc{\bbW}{{\mathbb W}}
 \nc{\bbX}{{\mathbb X}}
 \nc{\bbY}{{\mathbb Y}}
 \nc{\bbZ}{{\mathbb Z}}
 \nc{\bba}{{\mathbb a}}
 \nc{\bbb}{{\mathbb b}}
 \nc{\bbc}{{\mathbb c}}
 \nc{\bbd}{{\mathbb d}}
 \nc{\bbe}{{\mathbb e}}
 \nc{\bbf}{{\mathbb f}}
 \nc{\bbg}{{\mathbb g}}
 \nc{\bbh}{{\mathbb h}}
 \nc{\bbi}{{\mathbb i}}
 \nc{\bbk}{{\mathbb k}}
 \nc{\bbl}{{\mathbb l}}
 \nc{\bbm}{{\mathbb m}}
 \nc{\bbn}{{\mathbb n}}
 \nc{\bbo}{{\mathbb o}}
 \nc{\bbp}{{\mathbb p}}
 \nc{\bbq}{{\mathbb q}}
 \nc{\bbr}{{\mathbb r}}
 \nc{\bbs}{{\mathbb s}}
 \nc{\bbt}{{\mathbb t}}
 \nc{\bbu}{{\mathbb u}}
 \nc{\bbv}{{\mathbb v}}
 \nc{\bbw}{{\mathbb w}}
 \nc{\bbx}{{\mathbb x}}
 \nc{\bby}{{\mathbb y}}
 \nc{\bbz}{{\mathbb z}}

 \nc{\calA}{{\mathcal A}}
 \nc{\calB}{{\mathcal B}}
 \nc{\calC}{{\mathcal C}}
 \nc{\calD}{{\mathcal D}}
 \nc{\calE}{{\mathcal E}}
 \nc{\calF}{{\mathcal F}}
 \nc{\calG}{{\mathcal G}}
 \nc{\calH}{{\mathcal H}}
 \nc{\calI}{{\mathcal I}}
 \nc{\calJ}{{\mathcal J}}
 \nc{\tcalI}{{\tilde{\mathcal I}}}
 \nc{\tcalJ}{{\tilde{\mathcal J}}}
 \nc{\calK}{{\mathcal K}}
 \nc{\calL}{{\mathcal L}}
 \nc{\calM}{{\mathcal M}}
 \nc{\calN}{{\mathcal N}}
 \nc{\calO}{{\mathcal O}}
 \nc{\calP}{{\mathcal P}}
 \nc{\calQ}{{\mathcal Q}}
 \nc{\calR}{{\mathcal R}}
 \nc{\calS}{{\mathcal S}}
 \nc{\calT}{{\mathcal T}}
 \nc{\calU}{{\mathcal U}}
 \nc{\calV}{{\mathcal V}}
 \nc{\calW}{{\mathcal W}}
 \nc{\calX}{{\mathcal X}}
 \nc{\calY}{{\mathcal Y}}
 \nc{\calZ}{{\mathcal Z}}
  \nc{\cala}{{\mathcal a}}
 \nc{\calb}{{\mathcal b}}
 \nc{\calc}{{\mathcal c}}
 \nc{\cald}{{\mathcal d}}
 \nc{\cale}{{\mathcal e}}
 \nc{\calf}{{\mathcal f}}
 \nc{\calg}{{\mathcal g}}
 \nc{\calh}{{\mathcal h}}
 \nc{\cali}{{\mathcal i}}
 \nc{\calj}{{\mathcal j}}
 \nc{\calk}{{\mathcal k}}
 \nc{\call}{{\mathcal l}}
 \nc{\calm}{{\mathcal m}}
 \nc{\caln}{{\mathcal n}}
 \nc{\calo}{{\mathcal o}}
 \nc{\calp}{{\mathsf p}}
 \nc{\calq}{{\mathcal q}}
 \nc{\calr}{{\mathcal r}}
 \nc{\cals}{{\mathcal s}}
 \nc{\calt}{{\mathcal t}}
 \nc{\calu}{{\mathcal u}}
 \nc{\calv}{{\mathcal v}}
 \nc{\calw}{{\mathcal w}}
 \nc{\calx}{{\mathcal x}}
 \nc{\caly}{{\mathcal y}}
 \nc{\calz}{{\mathcal z}}

 \nc{\frakA}{{\mathfrak A}}
 \nc{\frakB}{{\mathfrak B}}
 \nc{\frakC}{{\mathfrak C}}
 \nc{\frakD}{{\mathfrak D}}
 \nc{\frakE}{{\mathfrak E}}
 \nc{\frakF}{{\mathfrak F}}
 \nc{\frakG}{{\mathfrak G}}
 \nc{\frakH}{{\mathfrak H}}
 \nc{\frakI}{{\mathfrak I}}
 \nc{\frakJ}{{\mathfrak J}}
 \nc{\frakK}{{\mathfrak K}}
 \nc{\frakL}{{\mathfrak L}}
 \nc{\frakM}{{\mathfrak M}}
 \nc{\frakN}{{\mathfrak N}}
 \nc{\frakO}{{\mathfrak O}}
 \nc{\frakP}{{\mathfrak P}}
 \nc{\frakQ}{{\mathfrak Q}}
 \nc{\frakR}{{\mathfrak R}}
 \nc{\frakS}{{\mathfrak S}}
 \nc{\frakT}{{\mathfrak T}}
 \nc{\frakU}{{\mathfrak U}}
 \nc{\frakV}{{\mathfrak V}}
 \nc{\frakW}{{\mathfrak W}}
 \nc{\frakX}{{\mathfrak X}}
 \nc{\frakY}{{\mathfrak Y}}
 \nc{\frakZ}{{\mathfrak Z}}
 \nc{\fraka}{{\mathfrak a}}
 \nc{\frakb}{{\mathfrak b}}
 \nc{\frakc}{{\mathfrak c}}
 \nc{\frakd}{{\mathfrak d}}
 \nc{\frake}{{\mathfrak e}}
 \nc{\frakf}{{\mathfrak f}}
 \nc{\frakg}{{\mathfrak g}}
 \nc{\frakh}{{\mathfrak h}}
 \nc{\fraki}{{\mathfrak i}}
 \nc{\frakj}{{\mathfrak j}}
 \nc{\frakk}{{\mathfrak k}}
 \nc{\frakl}{{\mathfrak l}}
 \nc{\frakm}{{\mathfrak m}}
 \nc{\frakn}{{\mathfrak n}}
 \nc{\frako}{{\mathfrak o}}
 \nc{\frakp}{{\mathfrak p}}
 \nc{\frakq}{{\mathfrak q}}
 \nc{\frakr}{{\mathfrak r}}
 \nc{\fraks}{{\mathfrak s}}
 \nc{\frakt}{{\mathfrak t}}
 \nc{\fraku}{{\mathfrak u}}
 \nc{\frakv}{{\mathfrak v}}
 \nc{\frakw}{{\mathfrak w}}
 \nc{\frakx}{{\mathfrak x}}
 \nc{\fraky}{{\mathfrak y}}
 \nc{\frakz}{{\mathfrak z}}

 \nc{\sha}{{\mbox{\cyr x}}}
\nc{\sofive}{{{\mathfrak {so}}(5)}}
\nc{\sofour}{{{\mathfrak {so}}(4)}}
\nc{\sothree}{{{\mathfrak {so}}(3)}}
\nc{\sosix}{{{\mathfrak {so}}(6)}}
\nc{\spsix}{{{\mathfrak {sp}}(6)}}
\nc{\soseven}{{{\mathfrak {so}}(7)}}
\nc{\sotwo}{{{\mathfrak {so}}(2)}}
\nc{\spone}{{{\mathfrak {sp}}(1)}}
\nc{\sptwo}{{{\mathfrak {sp}}(2)}}
\nc{\gtwo}{{\frg_2}}
\nc{\frg}{{\mathfrak g}}
 \nc{\slthree}{{{\mathfrak{sl}}(3)}}
 \nc{\slfour}{{{\mathfrak{sl}}(4)}}
 \nc{\slfive}{{{\mathfrak{sl}}(5)}}
 \nc{\sld}{{{\mathfrak{sl}}(d+1)}}
 \nc{\slr}{{{\mathfrak{sl}}(r+1)}}
 \nc{\slrr}{{{\mathfrak{sl}}(r+2)}}
 \nc{\gf}{{\varphi}}

 \nc{\uds}{{\underline{s}}}
\nc{\va}{{\vec a}}
\nc{\vb}{{\vec b}}
\nc{\vc}{{\vec c}}
\nc{\vdta}{{\vec \delta}}
\nc{\ve}{{\vec e}}
\nc{\vm}{{\vec m}}
\nc{\vp}{{\vec p}}
\nc{\vn}{{\vec n}}
\nc{\vmu}{{\vec \mu}}
\nc{\vr}{{\vec r}}
\nc{\vs}{{\vec s}}
\nc{\vt}{{\vec t}}
\nc{\vu}{{\vec u}}
\nc{\vx}{{\vec x}}
\nc{\vC}{{\vec C}}
\nc{\vv}{{\bf v}}
\nc{\cicc}[1]{{}_{{}^{ \bigcirc\hskip-1.2ex{#1}\hskip.3ex{}}}}
\nc{\cic}[1]{{}^{\bigcirc\hskip-1.25ex{#1}\hskip.25ex{}}}
\nc{\ccic}[1]{{}_{{}^{ \bigcirc\hskip-1.2ex{#1}\hskip.3ex{}}}}
\nc{\cci}[1]{{}_{{}^{ {\textstyle \bigcirc}\hskip-2.05ex{#1}\hskip-.35ex{}}}}
\nc{\ccicc}[1]{{}_{{}^{ {\textstyle \bigcirc}\hskip-1.55ex{#1}\hskip-0.1ex{}}}}
\nc{\x}{\rm{x}}

\begin{document}

\title[Witten volume formula]
{Witten volume formulas for semi-simple Lie algebras}

%\subjclass{Primary: 11R42,11F22; Secondary: 11S45,11M41,11Y16}
%\keywords{Zeta-function of Witten's type, reducibility.}

\author[Jianqiang {\sc ZHAO}]{{\sc Jianqiang} ZHAO}
\address{Department of Mathematics\\
Eckerd College\\
St. Petersburg, FL 33711\\  USA}
\email{zhaoj@eckerd.edu}
\urladdr{http://home.eckerd.edu/$\sim$zhaoj}

\maketitle

\begin{abstract}
In this paper we provide an algebraic derivation of the
explicit Witten volume formulas for
a few semi-simple Lie algebras by combining a combinatorial method
with the ideas used by Gunnells and Sczech in computation of
higher-dimensional Dedekind sums.
\end{abstract}

\bigskip

\section{Introduction}
In \cite{W} Witten related the volumes of the moduli spaces of
representations of the fundamental groups of two dimensional surfaces
to the special values of the following zeta function attached to
complex semisimple Lie algebras $\frg$ at positive integers:
$$\gz_W(s;\frg)=\sum_\gf \frac1{(\dim \gf)^s},$$
where $\gf$ runs over all finite dimensional irreducible representations of $\frg$.
By physics consideration Witten showed that for any positive integer $m$
$$\gz_W(2m;\frg)=c(2m;\frg) \pi^{2mr},$$
where $c(2m;\frg)\in\Q$ and $r$ is the number positive roots of $\frg$.
Such formulas are now called Witten volume formula.

The precise Witten volume formula
for $\slthree$ was obtained by Zagier \cite{Zag}
(and independently by Garoufalidis and Weinstein):
\begin{equation}\label{equ:sl3witten}
\gz_W(2m;\slthree)=\frac{4^{m+1}}{3} \sum_{\substack{0\le i\le 2m\\ i\equiv 0\text{ mod }2}}
{4m-i-1\choose 2m-1}\gz(i)\gz(6m-i).
\end{equation}
In \cite{GS}, Gunnells and Sczech studied higher-dimensional Dedekind sums and
established their reciprocity law. As one of the applications they could
derived the Witten volume formula for $\slfour$ precisely.

Matsumoto and his collaborators recently defined the multiple variable
analogs of $\gz_W(s;\frg)$ and studied some of their analytical and
arithmetical properties (see \cite{KMT1,KMT23,KMT4,MTs}):
$$
\gz_\frg(\{s_\ga\}_{\ga\in \gD_+}):= \sum_{m_1,\dots,m_\ell=1}^\infty
 \prod_{\ga\in \gD_+} \langle \ga^\vee, m_1\gl_1+\dots+m_\ell\gl_\ell\rangle^{-s_\ga}
$$
where for fixed set $\gD=\{\ga_1,\dots,\ga_\ell\}$ of fundamental roots
$\gD_+$ is set of all positive roots of $\frg$,
$\ga^\vee=2\ga/\langle\ga,\ga \rangle$ is the coroot,
and $\{\gl_1,\dots,\gl_\ell\}$ are the fundamental weights such that
$\langle \ga_i^\vee,\gl_j  \rangle=\gd_{i,j}$. By simple computation
$$
\gz_W(s;\frg)=M(\frg)^s  \gz_\frg( s,\dots,s),
\quad
\text{where} \quad
M(\frg)= \prod_{\ga\in \gD_+} \langle \ga^\vee,\gl_1+\dots+\gl_\ell\rangle.$$
With this multiple variable setup
Matsumoto et al. were able to obtain more general formulas
which include Witten volume formulas as special cases for
Lie algebras such as $\sofive$ and $\gtwo$. However, their computation
involves complicated analytical tools.

In this paper, we combine our combinatorial method developed in
\cite{Zgenzeta,Zso5,Z12th} and the technique of Gunnells and Sczech
to provide an algebraic proof of Witten volume formulas for $\sofive$ and $\gtwo$.

This paper is inspired by the work of Gunnells and Sczech \cite{GS}.
I want to thank them for their detailed explanation of
the part of their paper closely related to Witten zeta functions.

\section{The key ideas}
We briefly recall the setup in \cite[\S1]{GS}. Let $L$ be a lattice of rank $\ell\ge 1$
and $L^*=\text{Hom}_\Z(L,\Z)$. Denote by $\textbf{0}$ the zero linear form in $L^*$.
Let $r\ge \ell$. For $e=(e_1,\dots,e_r)\in \N^r$, $v\in L^*\otimes \R$ and
$\gs=(\gs_1,\dots,\gs_r)\in (L^*\setminus\{\textbf{0}\})^r$ Gunnells and Sczech
define the Dedekind sum
\begin{equation*}
D(L,\gs,e,v):=\frac1{(2\pi\sqrt{-1})^{\wt(e)}}\sum_{x\in L} \hskip-3.5ex {\phantom\sum}'
\frac{\exp(2\pi \sqrt{-1}\langle x,v\rangle)}{\langle x,\gs_1\rangle^{e_1} \cdots \langle x,\gs_r\rangle^{e_r}},
\end{equation*}
where $\langle\ ,\ \rangle:L\times L^*\to \Z$ is the pairing and $\sum'$ means
the terms with vanishing denominator are to be omitted.
When $L=\Z^\ell$ we represent each $v\in L^*\otimes \R$ by a vector in $\R^\ell$
so that $\langle (m_1,\dots,m_\ell),(v_1,\dots, v_\ell)\rangle=m_1v_1+\cdots+m_\ell v_\ell.$

Let $\ell$ be the rank of the semisimple Lie algebra $\frg$, $r=|\gD_+|$, and $W$ its Weyl group.
Define an $\ell\times r$ integral matrix $\gs(\frg)$ whose $j$-th column $v_j$ provides the coefficients
of $\ga_j\in \gD_+$ in terms of the fundamental roots in $\gD$. Let $e=(2m,\dots,2m)\in \N^r$.
Then by \cite[Prop.\ 8.4]{GS} we have
\begin{equation}\label{equ:gzD}
 \gz_W(2m,\frg)=(2\pi \sqrt{-1})^{2m r} \frac{M(\frg)^{2m}}{|W|}  D(\Z^\ell, \gs(\frg), e,\textbf{0}).
\end{equation}

In \cite{GS} Gunnells and Sczech demonstrated how one can use
the reciprocity law of higher-dimensional Dedekind sums to
derived the Witten volume formulas of some Lie algebras.
We are going to replace this tool by
the following simple combinatorial lemma (see \cite[p.~48]{Nie}).
\begin{lem} \label{lem:combLem}
Let $s,t$ be two positive integers. Let $x$ and $y$
be two non-zero real numbers such that $x+y\ne 0$. Then
$$ \frac{1}{x^sy^t} =
\sum_{a=0}^{s-1}  {t+a-1\choose a}\frac{1}{x^{s-a}(x+y)^{t+a}}
+\sum_{b=0}^{t-1} {s+b-1\choose b}\frac{1}{y^{t-b}(x+y)^{s+b}}.$$
\end{lem}

To demonstrate this idea we have the following key lemma to be used many times later
in the paper. Given any $\ell\times r$ matrix $\gs=(\gs_1,\dots,\gs_r)$ we denote
$\big((\gs_1)_{e_1},\dots,(\gs_r)_{e_r}\big)$ the new matrix obtained by repeating each
linear form $\gs_j$ exactly $e_j$ times, $j=1,\dots,r$.
For simplicity we further set
$$Z\big((\gs_1)_{e_1},\dots,(\gs_r)_{e_r}\big)=
(2 \pi  \sqrt{-1})^{e_1+\dots+e_r} D(\Z^\ell,(\gs_1,\dots,\gs_r),(e_1,\dots,e_r),\textbf{0}).$$
For example if
\begin{equation*}
M=\bigg(\begin{matrix}
1_{\phantom 4}    & 1_{\phantom 2}  \\
0_4  & 2_2
\end{matrix}\bigg)
=
\bigg(\begin{matrix}
1 & 1 & 1 & 1 & 1 & 1 \\
0 & 0 & 0 & 0 & 2 & 2
\end{matrix}\bigg)
\end{equation*}
then
$$Z(M) =-(2 \pi)^6 D\bigg(\Z^2,\bigg(\begin{matrix}
1 & 1  \\
0 & 2
\end{matrix}\bigg),(4,2),\textbf{0}\bigg).$$

\begin{lem} \label{lem:keyLem}
Suppose $a,b,c,d,e,f\in \Z$ such that $\gcd(a,b)=\gcd(c,d)=\gcd(e,f)=1$ and
$u\gl (a,b)+\gl (c,d)=(e,f)$ for some nonzero constant $\gl$ and
$|u|=1, 2$ or $1/2$. Let $\gd=\big|\begin{smallmatrix}
a & e \\
b & f
\end{smallmatrix} \big|.$ If $|\gd|=1,2$, and $|u\gd|=1,2$
then for all positive integers $i,j$ and $k$ such that $w=i+j+k$ is even we have
\begin{align}\label{equ:keylemma1}
Z\bigg(\begin{matrix}
a_{\phantom i} & c_{\phantom j} & e_{\phantom k} \\
b_i & d_j & f_k
\end{matrix}\bigg)=& (2\pi \sqrt{-1})^w \sum_{l=0}^i{i+j-l-1\choose j-1}
u^{i-l}\gl^{i+j-l}\frac{B_lB_{w-l}}{l!(w-l)!}\ga_{l,w}(\gd) \\
+& (2\pi \sqrt{-1})^w \sum_{l=0}^j {i+j-l-1\choose i-1}
u^{i} \gl^{i+j-l}\frac{B_lB_{w-l}}{l!(w-l)!} \ga_{l,w}(u\gd),  \label{equ:keylemma2}
\end{align}
where $\ga_{l,w}(\pm 1)=1$ and $\ga_{l,w}(\pm 2)=1-1/2^l-1/2^{w-l}+2/2^w.$
\end{lem}
\begin{proof} Clearly we have
\begin{equation}\label{equ:unimod}
\bigg|
\begin{matrix}
a & e \\
b & f
\end{matrix} \bigg|=\gd \qquad \Longrightarrow\qquad
\bigg|
\begin{matrix}
c & e \\
d & f
\end{matrix} \bigg|= -u\gd.
\end{equation}

For all integers $x$ and $y$ let $(x, y)^\bot=\{(m_1,m_2)\in \Z^2: x m_1+y m_2=0\}$.
By the combinatorial Lemma \ref{lem:combLem} we have
\begin{align}
Z
\bigg(\begin{matrix}
a_{\phantom i} & c_{\phantom j} & e_{\phantom k} \\
b_i & d_j & f_k
\end{matrix}\bigg)=& \sum_{l=1}^i{i+j-l-1\choose j-1}  u^{i-l}\gl^{i+j-l}
\left[Z
\bigg(\begin{matrix}
a_{\phantom l}  & e_{\phantom {w-l}}  \\
b_l & f_{w-l}
\end{matrix}\bigg) -
Z_{(c, d)^\bot} \bigg(\begin{matrix}
a_{\phantom l}  & e_{\phantom {w-l}}  \\
b_l & f_{w-l}
\end{matrix}\bigg) \right]  \notag \\
+& \sum_{l=1}^j {i+j-l-1\choose i-1} u^{i}\gl^{i+j-l}\left[Z
\bigg(\begin{matrix}
c_{\phantom l}  & e_{\phantom {w-l}} \\
d_l & f_{w-l}
\end{matrix}\bigg) -
Z_{(a, b)^\bot} \bigg(\begin{matrix}
c_{\phantom l}  & e_{\phantom {w-l}}  \\
d_l & f_{w-l}
\end{matrix}\bigg) \right]. \label{equ:sum=}
\end{align}
Here for any sub lattice $L$ of $\Z^2$ the sum $Z_L$ is the sum
$Z$ restricted to $L$. These terms in fact exactly correspond
to those appearing on the right hand of the reciprocity law \cite[(15)]{GS}.
The condition $\gcd(c,d)=1$ implies that
\begin{align*}
Z_{(c, d)^\bot}\bigg(\begin{matrix}
a_{\phantom l}  & e_{\phantom {w-l}}  \\
b_l & f_{w-l} \end{matrix} \bigg)
=&\sum_{\substack{m_1,m_2\in \Z,\ c m_1+d m_2=0\\
    am_1+bm_2\ne 0,\ em_1+fm_2\ne 0}} \frac1{(am_1+bm_2)^l(em_1+fm_2)^{w-l}}    \\
=&\sum_{N\in \Z^{{}^*},\ m_1=dN,\ m_2=-cN}\frac{u^l\gl^l}{(u\gl am_1+u\gl bm_2)^l(em_1+fm_2)^{w-l}} \\
=&\sum_{N\in \Z^{{}^*}}\frac{u^l\gl^l}{((e-\gl c)dN-(f-\gl d) cN)^l(edN-fcN)^{w-l}}
= \frac1{(u\gd)^w} \sum_{N\in \Z^{{}^*} }\frac{u^l\gl^l}{N^w}
\end{align*}
by \eqref{equ:unimod} since $w$ is even.
Hence by an easy binomial identity we get
\begin{equation}\label{equ:givel=0}
\sum_{l=0}^i{i+j-l-1\choose j-1} u^{i-l} \gl^{i+j-l}
Z_{(c, d)^\bot}\bigg(\begin{matrix}
a_{\phantom l}  & e_{\phantom {w-l}}  \\
b_l & f_{w-l}
\end{matrix}\bigg) ={i+j-1\choose j} \frac{u^i \gl^{i+j}}{(u\gd)^w} \sum_{N\in \Z^{{}^*}}\frac{1}{N^w}.
\end{equation}
Notice that for all integral matrix
$\gs={x \ y \choose s \ t}$ of determinant $\gd$
and positive integers $j,k$ with $j+k\ge 4$ we have  by \cite[(4)]{GS}
\begin{equation}\label{equ:detpm1}
Z \bigg(\begin{matrix}
x_{\phantom j}     & y_{\phantom {k}}   \\
s_j  & t_k
\end{matrix}\bigg)
= (2\pi \sqrt{-1})^{j+k}D\left(\Z^2,  \gs,
(j,k),\textbf{0}\right)=\frac{(2\pi \sqrt{-1})^{j+k}}{j!k!|\gd|}\sum_{z\in \Z^2/\gs \Z^2} \
     \mathcal{B}_j(\gs^{-1} z)  \mathcal{B}_k(\gs^{-1} z),
\end{equation}
where  $\mathcal{B}_j(x)$ are the Bernoulli polynomials.
When $u\gd=\pm 1$ the quantity in \eqref{equ:givel=0} provides exactly the $l=0$ term
in the sum of \eqref{equ:keylemma2} by the following
formula known to Euler: for even positive integer $w$
\begin{equation*}
 \sum_{N\in \Z^{{}^*}}\frac{1}{N^w} = 2\zeta(w)=-(2\pi)^w  \frac{B_{w}}{w!}.
\end{equation*}
When  $u\gd=\pm 2$ by the assumption $\gcd(c,d)=\gcd(e,f)=1$ and \eqref{equ:detpm1} we get
\begin{equation}\label{equ:gd=2}
Z\bigg(\begin{matrix}
c_{\phantom l}  & e_{\phantom {w-l}} \\
d_l & f_{w-l}
\end{matrix}\bigg)= \frac {(2\pi \sqrt{-1})^w} 2\frac{B_lB_{w-l}}{l!(w-l)!}\Big(1+\Big(\frac{2}{2^l}-1\Big)\Big(\frac{2}{2^{w-l}}-1\Big)\Big).
\end{equation}
When $l=0$ we find that \eqref{equ:gd=2} is equal to $B_{w}/(2^w w!)$ and
therefore \eqref{equ:givel=0} again provides exactly the $l=0$ term
in the sum of \eqref{equ:keylemma2}.

Similarly, the $l=0$ term in the sum of \eqref{equ:keylemma1} can be obtained by the sum of
$Z_{(a, b)^\bot}\big(\begin{smallmatrix}
c_{\phantom l}  & e_{\phantom {w-l}}  \\
d_l & f_{w-l}
\end{smallmatrix}\big)$ in \eqref{equ:sum=}.
This finishes the proof of the lemma.
\end{proof}

We can obtain \eqref{equ:sl3witten} immediately by applying the lemma to
$$Z\bigg(\begin{matrix}
1_{\phantom {2m}}     & 0_{\phantom {2m}}  & 1_{\phantom {2m}}  \\
0_{2m}  & 1_{2m} & 1_{2m}
\end{matrix}\bigg).$$

To aid our computation we represent the procedure in Lemma \ref{lem:keyLem}
by the following picture: the left is self-evident
while the right is more elegant by only recording the merged columns.
For example, in $1/xy(x+y)=1/x(x+y)^2+1/y(x+y)^2$ we may think $x$, $y$
and $x+y$ correspond to the first three columns and $1/x(x+y)^2$ is said
to be obtained by merge the 2nd column to the 3rd.
We can generalize Lemma \ref{lem:combLem} and apply it to any three linearly
dependent columns of a general matrix.
In the following picture, a circled column number between any two
sub-nodes signifies the column to which the two nodes are merged.
We call the tree in the right picture a \emph{computation tree}.
\begin{center}
\begin{tikzpicture}[scale=0.9]
\node (gs1) at  (0,2.7)  {$\bigg(\begin{matrix}
a_{\phantom i} & c_{\phantom j} & e_{\phantom k} \\
b_i & d_j & f_k
\end{matrix}\bigg)
$};
\node (B1) at (-1.5,0.8) {$\bigg(\begin{matrix}
a_{\phantom l} & e_{\phantom {w-l}} \\
b_l & f_{w-l}
\end{matrix}\bigg)$};
\node (B2) at (1.5,0.8) {$\bigg(\begin{matrix}
c_{\phantom l} & e_{\phantom {w-l}} \\
d_l & f_{w-l}
\end{matrix}\bigg)$};

\draw (gs1) to (B1);
\draw (gs1) to (B2);

\node (gs) at  (6.5,2.7)  {$\underset{\cic{1}\quad\ \cic{2}\quad\ \cic{3}}{\bigg(\begin{matrix}
a_{\phantom i} & c_{\phantom j} & e_{\phantom k} \\
b_i & d_j & f_k
\end{matrix}\bigg)}
$};
\node (A1) at (5,0.8) {2};
\node (A2) at (8,0.8) {1};
\draw (6.5,0.8) node {$\cicc{3}$};

\draw (gs) to (A1);
\draw (gs) to (A2);
\end{tikzpicture}
\end{center}

We have the following partial generalization of Lemma \ref{lem:keyLem}
whose proof is left to the interested reader.
\begin{lem} \label{lem:keyLem2}
Let $\gs=(\gs_1,\dots,\gs_r)$ be an $(r-1)\times r$ matrix with $r\ge 3$. Let
$e_1,\dots,e_r\in \N$ and put $s=e_1+e_2+e_3$. Suppose $u=\pm 1$,
$u\gl\gs_1+\gl\gs_2=\gs_3$ for some non-zero constant $\gl$ and
$\det(\gs_2,\dots,\gs_{r})=\pm 1$. If both $\gs_1$ and $\gs_2$ have at least
one component equal to $\pm 1$ then
\begin{align}\label{equ:keylemma3}
Z\big((\gs_1)_{e_1},\dots,(\gs_r)_{e_r}\big)
= (2\pi \sqrt{-1})^w \frac{B_{e_4}\dots B_{e_r}}{e_4!\dots e_r!} u^{e_1}\gl^{e_1+e_2}\Bigg\{&\sum_{l=0}^{e_1} {e_1+e_2-l-1\choose e_2-1}
\frac{u^{l}B_lB_{s-l}}{\gl^{l} l!(s-l)!} \\
+ \sum_{l=0}^{e_2}& {e_1+e_2-l-1\choose e_1-1}
\frac{B_lB_{s-l}}{\gl^{l}l!(s-l)!} \Bigg\}.\label{equ:keylemma4}
\end{align}
\end{lem}
\begin{rem}
It is incorrect to generalize Lemma \ref{lem:keyLem} to arbitrary rank
without the assumption that $\gs_1$ and $\gs_2$ have at least
one component equal to $\pm 1$.
\end{rem}

\begin{defn} If a node of binary tree has the following property
then we say it's a \emph{good parent}: either it is not a grandparent,
or every one of its descendants names
one of their children the same as its only sibling.
\end{defn}

\begin{prop} \label{prop:keyprop}
Let $\gs=(\gs_1,\dots,\gs_r)$ be an $\ell \times r$ matrix with $r\ge \ell+1 \ge 2$.
Suppose that every column has some component equal to $\pm 1$. Let
$e_1,\dots,e_r\in \N$ and put $s=e_1+e_2+e_3$. Assume $\gs$ is a good parent
in its computation binary tree.
Further assume that we have the top part of the
computation binary tree as follows if it has four grandchildren:

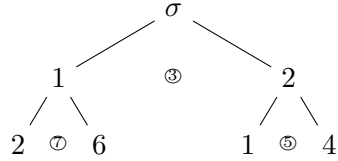
\begin{center}
\begin{tikzpicture}[scale=0.9]
\node (gs) at  (9.5,1.3)  {$\gs$};
\node (A1) at (7.8,0.3) {$1$};
\node (A2) at (11.2,0.3) {$2$};
\draw (9.5,0.3) node {$\cicc{3}$};
\node (B1) at (7.2,-0.7) {$2$};
\node (B2) at (8.4,-0.7) {$6$};
\draw (7.8,-0.7) node {$\cicc{7}$};
\node (B3) at (10.6,-0.7) {$1$};
\node (B4) at (11.8,-0.7) {$4$};
\draw (11.2,-0.7) node {$\cicc{5}$};

\draw (gs) to (A1);
\draw (A1) to (B1);
\draw (A1) to (B2);
\draw (gs) to (A2);
\draw (A2) to (B3);
\draw (A2) to (B4);
\end{tikzpicture}
\end{center}
Here, $3, 5$, and $7$ may or may not be the same but $5\ne 1,2,4$ and $7\ne 1,2,6.$
Suppose every node in the second to the last generation satisfies the conditions
in Lemma~\ref{lem:keyLem} (resp.\ Lemma~\ref{lem:keyLem2}) if it has rank two (resp.\ greater
than two). If $\gl_1\gs_1+\gl_2\gs_2=\gs_3$. Then
\begin{align*}
 Z\big((\gs_1)_{e_2},\dots,(\gs_r)_{e_r}\big)
=&\gl_1^{e_1}\gl_2^{e_2}\sum_{i=0}^{e_2} {e_1+e_2-i-1\choose e_1-1}
\frac1{\gl_2^i}Z\big((\gs_2)_{i},(\gs_3)_{s-i},(\gs_4)_{e_4},\dots,(\gs_r)_{e_r}\big) \\
+&\gl_1^{e_1}\gl_2^{e_2}\sum_{i=0}^{e_1} {e_1+e_2-i-1\choose e_2-1}
\frac1{\gl_1^i}Z\big((\gs_1)_{i},(\gs_3)_{s-i},(\gs_4)_{e_4},\dots,(\gs_r)_{e_r}\big) .
\end{align*}
\end{prop}

\begin{proof}
When $r=\ell+1$ the proposition follows from Lemma~\ref{lem:keyLem} and Lemma~\ref{lem:keyLem2})
by our assumption since now the second to the last generation is exactly $\gs$ itself.
Assume $r>\ell+1$. By Lemma \ref{lem:combLem} it's clear that
\begin{align*}
\ &Z\big((\gs_1)_{e_1},\dots,(\gs_r)_{e_r}\big)\\
=&\gl_1^{e_1}\gl_2^{e_2}\sum_{i=1}^{e_2} {e_1+e_2-i-1\choose e_1-1}
\frac1{\gl_2^i}\Big[Z\big((\gs_2)_{i},(\gs_3)_{s-i},(\gs_4)_{e_4},\dots,(\gs_r)_{e_r}\big) -Z_{\gs_1^\bot}\Big]\\
+&\gl_1^{e_1}\gl_2^{e_2}\sum_{i=1}^{e_1} {e_1+e_2-i-1\choose e_2-1}
\frac1{\gl_1^i}\Big[Z\big((\gs_1)_{i},(\gs_3)_{s-i},(\gs_4)_{e_4},\dots,(\gs_r)_{e_r}\big)-Z_{\gs_2^\bot}\Big] .
\end{align*}
With fixed $\ell$ we now use induction on $r$ to show that
\begin{align}
 \gl_2^{-i}Z_{\gs_1^\bot}\big((\gs_3)_{s},(\gs_4)_{e_4},\dots,(\gs_r)_{e_r}\big)
=&Z\big((\gs_1)_{i},(\gs_3)_{s-i},(\gs_4)_{e_4},\dots,(\gs_r)_{e_r}\big)\big|_{i=0}, \notag\\
 \gl_1^{-i} Z_{\gs_2^\bot}\big((\gs_3)_{s},(\gs_4)_{e_4},\dots,(\gs_r)_{e_r}\big)
=&Z\big((\gs_2)_{i},(\gs_3)_{s-i},(\gs_4)_{e_4},\dots,(\gs_r)_{e_r}\big)\big|_{i=0},\label{equ:ind2}
\end{align}
which yields the Lemma immediately by the simple combinatorial identities
\begin{equation*}
\sum_{i=1}^{e_1} {e_1+e_2-i-1\choose e_2-1} ={e_1+e_2-1\choose e_2},\quad
\sum_{i=1}^{e_2} {e_1+e_2-i-1\choose e_1-1} ={e_1+e_2-1\choose e_1}.
\end{equation*}
By the computation tree we may assume $\mu_1\gs_1+\mu_2\gs_4=\gs_5$.
Notice that for any $x\in \gs_1^\bot $ we have
\begin{align*}
 \gl_2\langle\gs_2,x\rangle=&\gl_1\langle\gs_1,x\rangle+\gl_2\langle\gs_2,x\rangle=\langle\gs_3,x\rangle,\\
  \mu_2\langle\gs_4,x\rangle=&\mu_1\langle\gs_1,x\rangle+\mu_2\langle\gs_4,x\rangle=\langle\gs_5,x\rangle.
H\end{align*}
ence (if $\cicc{3}=\cicc{5}$ then $e_5=0$)
\begin{align} \notag
\gl_2^{-i}Z_{\gs_1^\bot}\big((\gs_2)_{i},(\gs_3)_{s-i},(\gs_4)_{e_4},\dots,(\gs_r)_{e_r}\big)
=&Z_{\gs_1^\bot}\big((\gs_3)_{s},(\gs_4)_{e_4},\dots,(\gs_r)_{e_r}\big)\\
=&\mu_2^{e_4} Z_{\gs_1^\bot}Z\big((\gs_3)_{s},(\gs_5)_{e_4+e_5},(\gs_6)_{e_6},\dots,(\gs_r)_{e_r}\big).\label{equ:gs1bot}
\end{align}
On the other hand, by induction assumption we get
\begin{align*}
\ &Z\big((\gs_1)_{i},(\gs_3)_{s-i},(\gs_4)_{e_4},\dots,(\gs_r)_{e_r}\big) \\
=&\mu_1^i\mu_2^{e_4}\sum_{j=0}^{e_4} {e_4+i-j-1\choose i-1}
\frac1{\mu_2^j}Z\big((\gs_3)_{s-i},(\gs_4)_{j},(\gs_5)_{e_4+e_5+i-j},(\gs_6)_{e_6},\dots,(\gs_r)_{e_r}\big) \\
+& \mu_1^i\mu_2^{e_4} \sum_{j=0}^{i} {e_4+i-j-1\choose e_4-1}
\frac1{\mu_1^j}Z\big((\gs_1)_{j},(\gs_3)_{s-i},(\gs_5)_{e_4+e_5+i-j},(\gs_6)_{e_6},\dots,(\gs_r)_{e_r}\big).
\end{align*}
Taking $i=0$ in this expression we see that the first sum is vacuous because of the binomial coefficient
while the second sum is reduced to just one term:
\begin{multline*}
\mu_1^{i-j}\mu_2^{e_4} Z\big((\gs_1)_{j},(\gs_3)_{s-i},(\gs_5)_{e_4+e_5+i-j},(\gs_6)_{e_6},\dots,(\gs_r)_{e_r}\big)\big|_{i=j=0}\\
=\mu_2^{e_4} Z_{\gs_1^\bot}\big((\gs_3)_{s},(\gs_5)_{e_4+e_5},(\gs_6)_{e_6},\dots,(\gs_r)_{e_r}\big),
\end{multline*}
by induction assumption. Thus equation \eqref{equ:ind2} follows from \eqref{equ:gs1bot}.
The proof of \eqref{equ:ind2} is exactly the same. This concludes the proof of the proposition.
\end{proof}

\section{The $\sofive$ case}
Let $m\in \N$ and $n=2m$. By the above we can write
\begin{equation*}
\gz_\sofive(n,n,n,n)=\sum_{a,b=1}^\infty
\frac{1}{a^n b^n (a+b)^n (a+2b)^n}=
\frac{(2\pi)^{8m}}{8} D(\Z^2,\gs,(1,\dots,1),\textbf{0})
\end{equation*}
where $(1,\dots,1)\in \N^{4n}$, the matrix $\gs=\gs(n,n,n,n)$ and
\begin{equation*}
\gs(a,b,c,d)=\bigg(\begin{matrix}
1_{\phantom a}   & 0_{\phantom b}  & 1_{\phantom c}  & 1_{\phantom d}  \\
0_a & 1_b & 1_c & 2_d
\end{matrix}\bigg).
\end{equation*}

To prepare for the $\soseven$ case we first prove a
generalization of the $\sofive$ case by the following computation tree
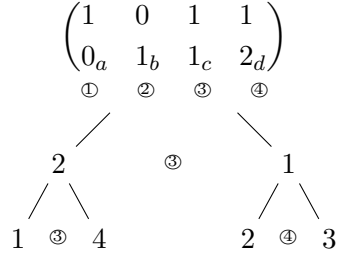
\begin{center}
\begin{tikzpicture}[scale=0.9]
\node (gs) at  (9.5,2.5)  {$\underset{\cic{1}\quad\ \cic{2}\quad\ \cic{3}\quad\ \cic{4}}{\bigg(\begin{matrix}
1_{\phantom a}   & 0_{\phantom b}  & 1_{\phantom c}  & 1_{\phantom d}  \\
0_a & 1_b & 1_c & 2_d
\end{matrix}\bigg)}
$};
\node (A1) at (7.8,0.8) {2};
\node (A2) at (11.2,0.8) {1};
\draw (9.5,0.8) node {$\cicc{3}$};
\node (B1) at (7.2,-0.3) {1};
\node (B2) at (8.4,-0.3) {4};
\draw (7.8,-0.3) node {$\cicc{3}$};
\node (B3) at (10.6,-0.3) {2};
\node (B4) at (11.8,-0.3) {3};
\draw (11.2,-0.3) node {$\cicc{4}$};

\draw (gs) to (A1);
\draw (A1) to (B1);
\draw (A1) to (B2);
\draw (gs) to (A2);
\draw (A2) to (B3);
\draw (A2) to (B4);
\end{tikzpicture}
\end{center}
\begin{thm}
Let $a,b,c,d\in \N$ and suppose at most one of them is $1$. Put
$\gb_{j,w}=B_jB_{w-j}/(j!(w-j)!).$ If $w=a+b+c+d$ is even then
\begin{align*}
Z\bigg(\begin{matrix}
1_{\phantom a} & 0_{\phantom a} & 1_{\phantom a}  & 1_{\phantom n}  \\
0_a  & 1_b & 1_c & 2_d
\end{matrix}\bigg)
=&\sum_{i=0}^b {a+b-i-1\choose a-1}\bigg(
\sum_{j=0}^i {w-d-j-1\choose w-d-i-1} \gb_{j,w}+
\sum_{j=0}^{w-d-i} {w-d-j-1\choose i-1} \gb_{j,w}\bigg)\\
+&\sum_{i=0}^a {a+b-i-1\choose b-1}\bigg(
\sum_{j=0}^d {d+i-j-1\choose i-1} \frac{\gb_{j,w}}{2^{d+i-j}}+
\sum_{j=0}^i {d+i-j-1\choose d-1} \frac{\gb_{j,w}}{2^{d+i-j}}\bigg)
\end{align*}
\end{thm}
\begin{proof} It's easy to check that all the $2\times 2$ minors of
$\bigg(\begin{matrix}
1_{\phantom a} & 0_{\phantom a} & 1_{\phantom a}  & 1_{\phantom n}  \\
0_a  & 1_b & 1_c & 2_d
\end{matrix}\bigg)$ have determinant $\pm 1$ or $\pm 2$. So we
can apply Proposition~\ref{prop:keyprop} and get
\begin{align*}
\ & Z\bigg(\begin{matrix}
1_{\phantom a} & 0_{\phantom a} & 1_{\phantom a}  & 1_{\phantom n}  \\
0_a  & 1_b & 1_c & 2_d
\end{matrix}\bigg)\\
=&\sum_{i=0}^b {a+b-i-1\choose a-1}
Z\bigg(\begin{matrix}
0_{\phantom a} & 1_{\phantom {w-d-i} }  & 1_{\phantom n}  \\
1_i & 1_{w-d-i}  & 2_d
\end{matrix}\bigg)
+\sum_{i=0}^a {a+b-i-1\choose b-1}
Z\bigg(\begin{matrix}
1_{\phantom a} & 1_{\phantom {w-d-i} }  & 1_{\phantom n}  \\
0_i  & 1_{w-d-i} & 2_d
\end{matrix}\bigg).
\end{align*}
The theorem now follows from Lemma \ref{lem:keyLem} directly.
\end{proof}

This implies the following as an immediate corollary. In the rest of the paper
for nonnegative integer $l$ and even integer $w\ge 4$ we
set $\gb_{l,w}= (2\pi)^w B_lB_{w-l}/(w-l)!/l!$.

\begin{thm}
Let $n$ be a positive even integer and $w=4n$. Set ${t \choose -1}=0$ for all $t$. Then
\begin{align*}
\frac{8 \gz_W(n,\sofive)}{6^n (2\pi)^{4n}}
=&\sum_{i=0}^n {2n-i-1\choose n-1}
\Bigg[\sum_{j=0}^{3n-i}   {3n-j-1\choose i-1}\beta_{j,w}
+  \sum_{j=0}^i   {3n-j-1\choose 3n-i-1}\beta_{j,w} \\
+& \sum_{j=0}^n \frac{1}{2^{n+i-j}} {n+i-j-1\choose i-1}\beta_{j,w}
+\sum_{j=0}^i \frac{1}{2^{n+i-j}} {n+i-j-1\choose n-1}\beta_{j,w} \Bigg].
\end{align*}
\end{thm}
\begin{proof} Clear.
\end{proof}
\begin{rem}
By exchange the order of summation in the theorem we see that
our formula agrees with that of Matsumoto et al. in \cite{KMT5}.
\end{rem}

\section{The $\gtwo$ case}
By definition we have $\gz_W(s;\gtwo)=120^s \gz_\gtwo(s,\dots,s)$ where
\begin{equation}\label{equ:g2defn}
 \gz_\gtwo(s_1,\dots,s_6)=\sum_{a,b=1}^\infty
 \frac{1}{a^{s_1}b^{s_2} (a+b)^{s_3} (a+2b)^{s_4}(a+3b)^{s_5}(2a+3b)^{s_6}}.
\end{equation}
In the rest of this section we fix a positive even integer $n=2m$. By \eqref{equ:gzD} we have
$$\gz_W(2m,\gtwo)=\frac{(2\pi)^{12m} 120^{2m} }{12}D(\Z^2, \gs, (1,\dots,1),\textbf{0})$$
where $(1,\dots,1)\in \N^{6n}$ and
\begin{equation*}
\underset{\ \phantom{\gs=} \cic{1}\quad\ \cic{2}\quad\ \cic{3}\quad\ \cic{4}\quad\ \cic{5}\quad\ \cic{6}}{ \gs=\Bigg(\begin{matrix}
1_{\phantom a} & 0_{\phantom b} & 1_{\phantom c} & 1_{\phantom d} & 1_{\phantom d} & 2_{\phantom d} \\
0_n & 1_n & 1_n & 2_n & 3_n & 3_n
\end{matrix}\Bigg).}
\end{equation*}
Similar to the case of $\sofive$ we can proceed using the following picture:

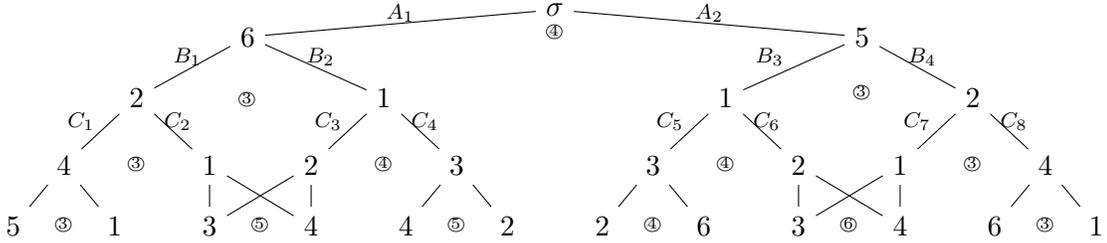
\begin{center}
\begin{tikzpicture}[scale=0.9]
\node (gs) at (0,0.8) {$\gs$};
\node (A1) at (-4.5375,0.4) {6};
\node (A2) at (4.5375,0.4)  {5};
\draw (0,0.45)   node  {$\cicc{4}$};
\node (B1) at (-6.175,-0.5) {2};
\node (B2) at (-2.525,-0.5) {1};
\draw (-4.5375,-0.55)   node {$\cicc{3}$};
\node (B3) at (2.525,-0.5) {1};
\node (B4) at (6.175,-0.5) {2};
\draw (4.5375,-0.45)   node  {$\cicc{3}$};

\node (C1) at (-7.25,-1.5) {4};
\node (C2) at (-5.1,-1.5) {1};
\draw (-6.175,-1.5)  node  {$\cicc{3}$};
\node (D1) at (-8,-2.4) {5};
\node (D2) at (-6.5,-2.4) {1};
\draw (-7.25,-2.4)  node  {$\cicc{3}$};
\node (D3) at (-5.1,-2.4) {3};
\node (D4) at (-3.6,-2.4) {4};
\draw (-4.35,-2.4)  node  {$\cicc{5}$};

\node (C3) at (-3.6,-1.5) {2};
\node (C4) at (-1.45,-1.5) {3};
\draw (-2.525,-1.5)  node {$\cicc{4}$};
\node (D5) at (-2.2,-2.4) {4};
\node (D6) at (-0.7,-2.4) {2};
\draw (-1.45,-2.4)  node  {$\cicc{5}$};

\node (C5) at (1.45,-1.5) {3};
\node (C6) at (3.6,-1.5) {2};
\draw (2.525,-1.5) node {$\cicc{4}$};
\node (D7) at (0.7,-2.4) {2};
\node (D8) at (2.2,-2.4) {6};
\draw (1.45,-2.4) node {$\cicc{4}$};

\node (C7) at (5.1,-1.5) {1};
\node (C8) at (7.25,-1.5) {4};
\draw (6.175,-1.5) node {$\cicc{3}$};
\node (D9) at (3.6,-2.4) {3};
\node (D10) at (5.1,-2.4) {4};
\draw (4.35,-2.4) node {$\cicc{6}$};

\node (D11) at (6.5,-2.4) {6};
\node (D12) at (8,-2.4) {1};
\draw (7.25,-2.4) node {$\cicc{3}$};

\draw (gs) to (A1) node[midway] {${\scriptstyle A_1}\atop\ $};
\draw (gs) to (A2) node[midway] {${\scriptstyle A_2}\atop\ $};
\draw (A1) to (B1) node[midway] {${{\scriptstyle B_1}\atop }\ $};
\draw (A1) to (B2) node[midway] {$\ {\scriptstyle B_2}\atop $};
\draw (A2) to (B3) node[midway] {${ {\scriptstyle B_3}\atop}\quad\ \   $};
\draw (A2) to (B4) node[midway] {$\ \scriptstyle B_4 \atop \ $};
\draw (B1) to (C1) node[midway] {${{\scriptstyle C_1}\atop }\quad\  $};
\draw (B1) to (C2) node[midway] {$\ {\scriptstyle C_2}\atop $};
\draw (B2) to (C3) node[midway] {${\scriptstyle C_3\atop \ }\quad\  $};
\draw (B2) to (C4) node[midway] {$\ {\scriptstyle C_4 \atop \ }$};
\draw (B3) to (C5) node[midway] {${\scriptstyle C_5\atop \ }\quad\ $};
\draw (B3) to (C6) node[midway] {$\ {\scriptstyle C_6 \atop \ } $};
\draw (B4) to (C7) node[midway] {${\scriptstyle C_7\atop \ }\quad\ $};
\draw (B4) to (C8) node[midway] {$\ {\scriptstyle C_8 \atop \ } $};
\draw (C1) to (D1);
\draw (C1) to (D2);
\draw (C2) to (D3);
\draw (C2) to (D4);
\draw (C3) to (D3);
\draw (C3) to (D4);
\draw (C4) to (D5);
\draw (C4) to (D6);
\draw (C5) to (D7);
\draw (C5) to (D8);
\draw (C6) to (D9);
\draw (C6) to (D10);
\draw (C7) to (D9);
\draw (C7) to (D10);
\draw (C8) to (D11);
\draw (C8) to (D12);
\end{tikzpicture}
\end{center}

Set $\gd_2(l)=1+2^{1-6n}-2^{-l}-2^{l-6n}$, $\gb_{l,w}= (2\pi)^w B_lB_{w-l}/(w-l)!/l!$ and
$\gb_{l,w}''=  \gd_2(l) \gb_{l,w}.$ We thus get
\begin{equation}\label{equ:A1A2}
Z\bigg(\begin{matrix}
1_{\phantom a} & 0_{\phantom b} & 1_{\phantom c} & 1_{\phantom d} & 1_{\phantom d} & 2_{\phantom d} \\
0_n & 1_n & 1_n & 2_n & 3_n & 3_n
\end{matrix}\bigg)=
\sum_{i=0}^n{2n-i-1\choose n-1}\frac{A_1(i)+A_2(i)}{3^{2n-i}},
\end{equation}
where
$$A_1(i)=Z\bigg(\begin{matrix}
1_{\phantom i}  & 0_{\phantom b} & 1_{\phantom c} & 1_{\phantom {3n-i}} & 1_{\phantom d}   \\
0_n & 1_n & 1_n & 2_{3n-i} & 3_i
\end{matrix}\bigg), \quad A_1(i)=Z\bigg(\begin{matrix}
1_{\phantom i}  & 0_{\phantom b} & 1_{\phantom c} & 1_{\phantom {3n-i}} & 1_{\phantom d}   \\
0_n & 1_n & 1_n & 2_{3n-i} & 3_i
\end{matrix}\bigg).$$
Note that $\gs$ is not a good parent so we have to check the two perpendicular terms
can indeed be absorbed into the summand when setting $i=0$ in \eqref{equ:A1A2}.
This is not too difficult
after finding out the explicit expressions of $A_1$ and $A_2$ both of which are
good parents. Therefore it follows from Proposition~\ref{prop:keyprop} that
\begin{align*}
A_1(i)=&\sum_{j=0}^n {2n-j-1\choose n-1}\Big(B_1(i,j)+B_2(i,j)\Big),\\
A_2(i)=&\sum_{j=0}^n {2n-j-1\choose n-1}\Big(B_3(i,j)+B_4(i,j)\Big),
\end{align*}
where
\begin{align*}
B_1(i,j)=Z\bigg(\begin{matrix}
1_{\phantom i} & 1_{\phantom {3n-j}} & 1_{\phantom {3n-i}} & 1_{\phantom d}   \\
0_j & 1_{3n-j} & 2_{3n-i} & 3_i
\end{matrix}\bigg)=&
\sum_{k=0}^j{3n-i+j-k-1\choose 3n-i-1} \frac{C_1(i,k)}{2^{3n-i+j-k}}\\
&+\sum_{k=0}^{3n-i}{3n-i+j-k-1\choose j-1}\frac{C_2(i,k)}{2^{3n-i+j-k}},\\
B_2(i,j)
=Z\bigg(\begin{matrix}
0_{\phantom i} & 1_{\phantom {3n-j}} & 1_{\phantom {3n-i}} & 1_{\phantom d}   \\
1_j & 1_{3n-j} & 2_{3n-i} & 3_i
\end{matrix}\bigg)=&
\sum_{k=0}^{3n-j}{3n-k-1\choose j-1} C_3(i,k)+
\sum_{k=0}^j {3n-k-1\choose 3n-j-1}C_4(i,k) ,\\
B_3(i,j)
=Z\bigg(\begin{matrix}
0_{\phantom i}  & 1_{\phantom {3n-j}} & 1_{\phantom {3n-i}} & 2_{\phantom d}   \\
1_j  & 1_{3n-j} & 2_{3n-i} & 3_i
\end{matrix}\bigg)=&
\sum_{k=0}^j {3n-k-1\choose 3n-j-1}  C_3(i,k)
+\sum_{k=0}^{3n-j}{3n-k-1\choose j-1} C_4(i,k) ,\\
B_4(i,j)
=Z\bigg(\begin{matrix}
1_{\phantom i}  & 1_{\phantom {3n-j}} & 1_{\phantom {3n-i}} & 2_{\phantom d}   \\
0_j  & 1_{3n-j} & 2_{3n-i} & 3_i
\end{matrix}\bigg)=&\sum_{k=0}^{3n-i}{3n-i+j-k-1\choose j-1}\frac{C_4(i,k)}{2^{3n-i+j-k}}\\
+&\sum_{k=0}^j{3n-i+j-k-1\choose 3n-i-1} \frac{C_8(i,k)}{2^{3n-i+j-k}},
\end{align*}
where, by using Lemma \ref{lem:keyLem} we have
{\allowdisplaybreaks
\begin{align*}
C_1(i,k)=Z
\bigg(\begin{matrix}
1_{\phantom i} & 1_{\phantom {6n-k-i}} & 1_{\phantom i}   \\
0_k & 1_{6n-k-i} & 3_i
\end{matrix}\bigg)=&
\sum_{l=0}^{i}{i+k-l-1\choose k-1}\frac{2^k\gd_2(l) }{3^{i+k-l}}
 \gb_{l,w}
+\sum_{l=0}^{k}{i+k-l-1\choose i-1}\frac{2^{k-l}}{3^{i+k-l}}
\gb_{l,w},\\
C_2(i,k)=Z
\bigg(\begin{matrix}
1_{\phantom {6n-k-i}} & 1_{\phantom k} & 1_{\phantom i}  \\
1_{6n-k-i} & 2_k & 3_i
\end{matrix}\bigg)=&
\sum_{l=0}^i {i+k-l-1\choose k-1} (-1)^i 2^{k} \gd_2(l)\gb_{l,w}
+\sum_{l=0}^k {i+k-l-1\choose i-1} (-1)^i 2^{k-l} \gb_{l,w} ,\\
C_3(i,k)=Z
\bigg(\begin{matrix}
1_{\phantom k} & 1_{\phantom {6n-k-i}} & 1_{\phantom i}  \\
1_k & 2_{6n-k-i} & 3_i
\end{matrix}\bigg)=&
\sum_{l=0}^i {i+k-l-1\choose k-1} (-1)^i \frac{\gb_{l,w}}{2^{k+i-l}}
+\sum_{l=0}^k {i+k-l-1\choose i-1} (-1)^i  \frac{\gb_{l,w}}{2^{k+i-l}} ,\\
C_4(i,k)=Z
\bigg(\begin{matrix}
0_{\phantom i} & 1_{\phantom {6n-k-i}} & 1_{\phantom d}   \\
1_k & 2_{6n-k-i} & 3_i
\end{matrix}\bigg)=&
\sum_{l=0}^{k}{i+k-l-1\choose i-1} (-1)^k \gb_{l,w}
+\sum_{l=0}^{i}{i+k-l-1\choose k-1}(-1)^k  \gb_{l,w} ,\\
C_8(i,k)=Z
\bigg(\begin{matrix}
1_{\phantom i} & 1_{\phantom {6n-k-i}} & 2_{\phantom i}   \\
0_k & 1_{6n-k-i} & 3_i
\end{matrix}\bigg)=&
\sum_{l=0}^{k}{i+k-l-1\choose i-1}  \frac{\gb_{l,w}}{3^{k+i-l}}
+\sum_{l=0}^{i}{i+k-l-1\choose k-1} \frac{\gb_{l,w}}{3^{k+i-l}}.
\end{align*}}

Putting everything together we finally arrive at
\begin{thm}
Let $n$ be a positive even integer. Let $w=6n$ and set ${t \choose -1}=0$ for all $t$. Write
$\gb_{l,w}=(2\pi)^w B_lB_{w-l}/(w-l)!/l!$. Then
{\allowdisplaybreaks
\begin{align*}
\frac{12\gz_W(n,\gtwo)}{120^n}= &
\sum_{i=0}^n{2n-i-1\choose n-1} 3^{i-2n} \sum_{j=0}^n{2n-j-1\choose n-1}\times
 \\
\times\bigg\{& \sum_{k=0}^j{3n-i+j-k-1\choose 3n-i-1}  2^{i+k-3n-j}
\sum_{t=1}^2  \sum_{l=0}^{\gl_t} {i+k-l-1\choose i+k-\gl_t-1}
    \frac{1+2^k\gd_t(l)}{3^{i+k-l}} \gb_{l,w}
\\
+& \sum_{k=0}^{3n-i}{3n-i+j-k-1\choose j-1}
\sum_{t=1}^2  \sum_{l=0}^{\gl_t} {i+k-l-1\choose i+k-\gl-1}\big((-1)^k+ (-1)^i 2^{k}\gd_t(l)\big)\gb_{l,w}
\\
+& \sum_{k=0}^{3n-j}{3n-k-1\choose j-1}
 \sum_{\gl=i,k}  \sum_{l=0}^{\gl} {i+k-l-1\choose i+k-\gl-1}\big((-1)^k+ (-1)^i2^{l-i-k}\big)\gb_{l,w}
\\
+& \sum_{k=0}^j {3n-k-1\choose 3n-j-1}
  \sum_{\gl=i,k}  \sum_{l=0}^{\gl} {i+k-l-1\choose i+k-\gl-1}\big((-1)^k+ (-1)^i2^{l-i-k}\big)\gb_{l,w}
 \bigg\},
\end{align*}}
where $\gl_1=k$, $\gd_1(l)=2^{-l}$, $\gl_2=i$, and
$\gd_2(l)=1+2^{1-6n}-2^{-l}-2^{l-6n}$.
\end{thm}
\begin{rem}
Although we can not verify the agreement of our theorem with \cite{KMT4} we are sure
their result will follow by choosing another computation tree. However, we find
our data for $\gz_W(2m,\gtwo)$ ($m\le 10$) agree with those in \cite{KMT4}.
We have also verified numerically the correctness of these values by using
the definition \eqref{equ:g2defn}.
\end{rem}

\section{The cases of $\soseven$ and $\spsix$}
By definition
\begin{equation}
\label{equ:definso7}
   \gz_W(n;\soseven)=720^n  \sum_{m_1,m_2,m_3=1}^\infty
   \left(\frac{ 1/\big(m_1 m_2 m_3 (m_1+m_2) (m_2+m_3) (2m_2+m_3)\big)}
     {(2m_1+2m_2+m_3) (m_1+2m_2+m_3) (m_1+m_2+m_3) } \right)^{n},
\end{equation}
and
\begin{equation*}
   \gz_W(n;\spsix)=720^n
   \sum_{m_1,m_2,m_3=1}^\infty
   \left(\frac{ 1/\big(m_1 m_2 m_3(m_1+m_2)(m_2+m_3)(m_2+2m_3)\big)}
  {(m_1+m_2+m_3) (m_1+m_2+2m_3) (m_1+2m_2+2m_3) }\right)^{n}.
\end{equation*}
The corresponding matrix to $\soseven$ is
\begin{equation*}
\underset{\ \phantom{\gs=} \cic{1}\quad\ \cic{2}\quad\ \cic{3}\quad\ \cic{4}\quad\ \cic{5}\quad\ \cic{6}\quad\ \cic{7}\quad\ \cic{8}\quad\ \cic{9} }{ \gs=\Bigg(\begin{matrix}
1_{\phantom a} & 0_{\phantom b} & 0_{\phantom b} & 1_{\phantom c} &0_{\phantom b} & 0_{\phantom b} & 2_{\phantom d} & 1_{\phantom d} & 1_{\phantom d} \\
0_{\phantom n} & 1_{\phantom n} & 0_{\phantom n} & 1_{\phantom n} & 1_{\phantom n} & 2_{\phantom n} & 2_{\phantom n} & 2_{\phantom n} & 1_{\phantom n} \\
0_n & 0_n & 1_n & 0_n & 1_n & 1_n & 1_n & 1_n & 1_n
\end{matrix}\Bigg).}
\end{equation*}
Given four column vectors $\gs_1,\dots,\gs_4$ let $S(\gs_1,\dots,\gs_4)$ be the set of
the four possible choices of three columns. Then every triple columns of the following
are linearly dependent: \begin{equation}\label{equ:lincol}
\{(1,5,9),(4,5,8)\}\cup  S(1,6,7,8)\cup  S(2,3,5,6)\cup  S(3,4,7,9).
\end{equation}
These fourteen dependencies are the only 3-column dependencies and will be
used critically in the computation tree of $\gz_W(n,\soseven)$ for even positive integer $n$.
If we consider the matrix corresponding to $\spsix$
\begin{equation*}
\underset{  \cic{1}\quad\ \cic{2}\quad\ \cic{3}\quad\ \cic{4}\quad\ \cic{5}\quad\ \cic{6}\quad\ \cic{7}\quad\ \cic{8}\quad\ \cic{9} }{  \Bigg(\begin{matrix}
1_{\phantom a} & 0_{\phantom b} & 0_{\phantom b} & 1_{\phantom c} &0_{\phantom b} & 0_{\phantom b} & 1_{\phantom d} & 1_{\phantom d} & 1_{\phantom d} \\
0_{\phantom n} & 1_{\phantom n} & 0_{\phantom n} & 1_{\phantom n} & 1_{\phantom n} & 1_{\phantom n} & 2_{\phantom n} & 1_{\phantom n} & 1_{\phantom n} \\
0_n & 0_n & 1_n & 0_n & 1_n & 2_n & 2_n & 2_n & 1_n
\end{matrix}\Bigg)}
\end{equation*}
we find the following similar combinations of columns:
$$\{(2,7,8), (4,6,7)\}\cup  S(2,3,5,6)\cup  S(1,5,7,9)\cup  S(3,4,8,9).$$
But these can be obtained exactly from \eqref{equ:lincol} by the permutation $(12)(398657)$.

We may construct the following computation tree for $\soseven$ as in the previous cases.
A similar computation tree for $\spsix$ can then be obtained by applying the permutation $(12)(398657)$.
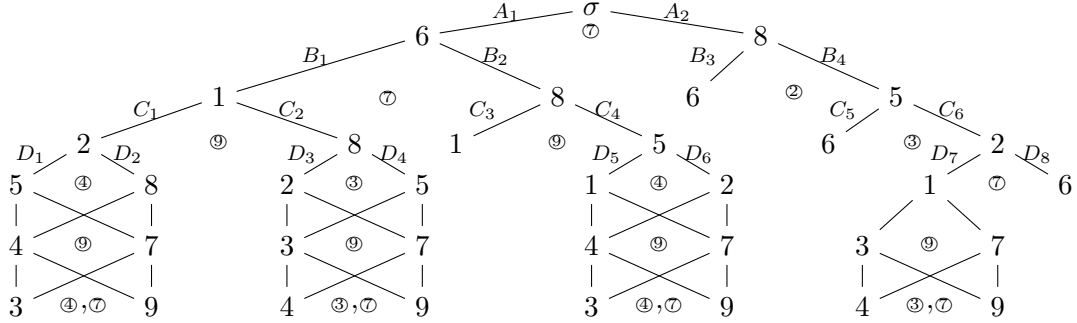
\begin{center}
\begin{tikzpicture}[scale=0.9]
\node (gs) at (0.5,0.8) {$\gs$};
\node (A1) at (-2,0.4) {6};
\node (A2) at (3,0.4)  {8};
\draw (0.5,0.45)   node  {$\cicc{7}$};
\node (B1) at (-5,-0.5) {1};
\node (B2) at (0,-0.5) {8};
\draw (-2.5,-0.55)   node {$\cicc{7}$};
\node (B3) at (2,-0.5) {6};
\node (B4) at (5,-0.5) {5};
\draw (3.5,-0.45)   node  {$\cicc{2}$};
\node (C1) at (-7,-1.2) {2};
\node (C2) at (-3,-1.2) {8};
\draw (-5,-1.2)  node  {$\cicc{9}$};
\node (C3) at (-1.5,-1.2) {1};
\node (C4) at (1.5,-1.2) {5};
\draw (0,-1.2)  node {$\cicc{9}$};
\node (C6) at (6.5,-1.2) {2};
\node (C5) at (4,-1.2) {6};
\draw (5.25,-1.2) node {$\cicc{3}$};
\node (D1) at (-8,-1.8) {5};
\node (D2) at (-6,-1.8) {8};
\draw (-7,-1.8) node {$\cicc{4}$};
\node (D3) at (-4,-1.8) {2};
\node (D4) at (-2,-1.8) {5};
\draw (-3,-1.8) node {$\cicc{3}$};
\node (D5) at (0.5,-1.8) {1};
\node (D6) at (2.5,-1.8) {2};
\draw (1.5,-1.8) node {$\cicc{4}$};
\node (D7) at (5.5,-1.8) {1};
\node (D8) at (7.5,-1.8) {6};
\draw (6.5,-1.8) node {$\cicc{7}$};

\node (E1) at (-8,-2.7) {4};
\node (E2) at (-6,-2.7) {7};
\draw (-7,-2.7) node {$\cicc{9}$};
\node (E3) at (-4,-2.7) {3};
\node (E4) at (-2,-2.7) {7};
\draw (-3,-2.7) node {$\cicc{9}$};
\node (E5) at (0.5,-2.7) {4};
\node (E6) at (2.5,-2.7) {7};
\draw (1.5,-2.7) node {$\cicc{9}$};
\node (E7) at (4.5,-2.7) {3};
\node (E8) at (6.5,-2.7) {7};
\draw (5.5,-2.7) node {$\cicc{9}$};

\node (F1) at (-8,-3.6) {3};
\node (F2) at (-6,-3.6) {9};
\draw (-7,-3.6) node {$\cicc{4}$,$\cicc{7}$};
\node (F3) at (-4,-3.6) {4};
\node (F4) at (-2,-3.6) {9};
\draw (-3,-3.6) node {$\cicc{3}$,$\cicc{7}$};
\node (F5) at (0.5,-3.6) {3};
\node (F6) at (2.5,-3.6) {9};
\draw (1.5,-3.6) node {$\cicc{4}$,$\cicc{7}$};
\node (F7) at (4.5,-3.6) {4};
\node (F8) at (6.5,-3.6) {9};
\draw (5.5,-3.6) node {$\cicc{3}$,$\cicc{7}$};

\draw (gs) to (A1) node[midway] {${\scriptstyle A_1}\atop\ $};
\draw (gs) to (A2) node[midway] {${\scriptstyle A_2}\atop\ $};
\draw (A1) to (B1) node[midway] {${{\scriptstyle B_1}\atop }\ $};
\draw (A1) to (B2) node[midway] {$\ {\scriptstyle B_2}\atop $};
\draw (A2) to (B3) node[midway] {${ {\scriptstyle B_3}\atop}\quad\ \   $};
\draw (A2) to (B4) node[midway] {$\ \scriptstyle B_4 \atop \ $};
\draw (B1) to (C1) node[midway] {${{\scriptstyle C_1}\atop }\ $};
\draw (B1) to (C2) node[midway] {$\ {\scriptstyle C_2}\atop $};
\draw (B2) to (C3) node[midway] {${\scriptstyle C_3 \atop \ }\quad\ \ $};
\draw (B2) to (C4) node[midway] {${\scriptstyle C_4\atop \ }$};
\draw (B4) to (C6) node[midway] {$ \ {\scriptstyle C_6 \atop \ } $};
\draw (B4) to (C5) node[midway] {${\scriptstyle C_5\atop \ }\quad\ $};
\draw (C1) to (D1) node[midway] {${\scriptstyle D_1 \atop \ }\quad\ $};
\draw (C1) to (D2) node[midway] {$\ \ {\scriptstyle D_2\atop \ }$};
\draw (C2) to (D3) node[midway] {${\scriptstyle D_3\atop \ }\quad\ $};
\draw (C2) to (D4) node[midway] {$\ {\scriptstyle D_4 \atop \ } $};
\draw (C4) to (D5) node[midway] {${\scriptstyle D_5\atop \ }\quad\ $};
\draw (C4) to (D6) node[midway] {$\ {\scriptstyle D_6 \atop \ } $};
\draw (C6) to (D7) node[midway] {${\scriptstyle D_7\atop \ }\quad\ $};
\draw (C6) to (D8) node[midway] {$\ {\scriptstyle D_8 \atop \ } $};
\draw (E1) to (D1);
\draw (E1) to (D2);
\draw (E2) to (D1);
\draw (E2) to (D2);
\draw (E3) to (D3);
\draw (E3) to (D4);
\draw (E4) to (D3);
\draw (E4) to (D4);
\draw (E5) to (D5);
\draw (E5) to (D6);
\draw (E6) to (D5);
\draw (E6) to (D6);
\draw (E7) to (D7);
\draw (E8) to (D7);
\draw (E1) to (F1);
\draw (E1) to (F2);
\draw (E2) to (F1);
\draw (E2) to (F2);
\draw (E3) to (F3);
\draw (E3) to (F4);
\draw (E4) to (F3);
\draw (E4) to (F4);
\draw (E5) to (F5);
\draw (E5) to (F6);
\draw (E6) to (F5);
\draw (E6) to (F6);
\draw (E7) to (F7);
\draw (E8) to (F7);
\draw (E7) to (F8);
\draw (E8) to (F8);
\end{tikzpicture}
\end{center}
Note that the column to which two sub-nodes are merged to may not be unique.
But it should not affect the final result. For example, at the very
beginning we can merge the 6th and 8th column to either the 1st
or 7th column. We will choose the 7th column in the computation tree.
On the other hand, if we go down the path $8-5-2-1-3$
then the 4th and 9th column can only be collapsed to the 7th column
since the 3rd was removed already.

It turns out that even though $\soseven$ and $\sosix$ are not isomorphic Lie algebras,
the computation of $\gz_W(n,\spsix)$ is almost exactly the same as that
of $\gz_W(n,\soseven)$. By following the above computation tree we have

\begin{thm} Let $n$ be an even positive integer. Set ${t \choose -1}=0$ for all $t$. Then
\begin{align*}
\frac{48\gz_W(n,\soseven)}{(-720)^n (2\pi)^{9n}}
=\frac{48\gz_W(n,\spsix)}{(-720)^n (2\pi)^{9n}}
=& 2^n\sum_{i=0}^n{2n-i-1\choose n-1}\Bigg\{ \sum_{j=0}^i {n+i-j-1\choose n-1}
    \bigg(\frac{B(j)}{2^i}+ \frac{B''(i,j)}{(-1)^i}\bigg)\\
+&  \sum_{j=0}^n {n+i-j-1\choose i-1} \bigg(\frac{B'(j,n,j)}{2^i}+\frac{B'(i+n,j,j)}{(-1)^{i+j}}\bigg)\Bigg\},\\
\end{align*}
where by writing $C'\{i,k\}=C'(n+i,3n+i-k,3n+i-k,2n+i-k)$
\begin{align*}
B(j)=&
\sum_{k=0}^j {n+j-k-1\choose n-1}C(j,j-k,k,j) +
\sum_{k=0}^n {n+j-k-1\choose j-1} \frac{C(j,j-k,n,j)}{(-1)^k},\\
B'(a,b,j)=&
\sum_{k=0}^b {n+j-k-1\choose n+j-b-1} C(a,j-k,k,j)+
\sum_{k=0}^{n+j-b} {n+j-k-1\choose b-1} C'(a,k,n,j),\\
B''(i,j)=&
\sum_{k=0}^{2n+i-j} {2n+i-k-1\choose j-1} \frac{2^{2n+i-j-k}C'\{i,k\}}{(-1)^{i-j-k}}
+
\sum_{k=0}^j {2n+i-k-1\choose 2n+i-j-1} \frac{2^{2n+i-j}C''(i,k)}{(-1)^{i-j} },
\end{align*}
and by setting $\mu=n+a-j$ and $\nu=n+b-c$
\begin{align*}
C(a,b,c,j)=&
\sum_{l=0}^c {n+a-b-l-1\choose n+a-b-c-1} \frac{D(a,b,n,j,l)}{(-1)^{n+a-b-c}}
+\sum_{l=0}^{n+a-b-c} {n+a-b-l-1\choose c-1} \frac{D(a,b,n,j,l)}{(-1)^{n+a-b-c-l}},\\
C'(a,b,c,j)=&
\sum_{l=0}^{\mu} {\mu+\nu-l-1\choose \nu-1} D(a,b,c,j,l)
+\sum_{l=0}^{\nu} {\mu+\nu-l-1\choose \mu-1} D(a,b,c,j,l),\\
C''(i,k)=&\sum_{l=0}^k {n+k-l-1\choose n-1} 2^n D'(i,k,l)
+\sum_{l=0}^n {n+k-l-1\choose k-1} 2^{n-l} D(i+l-k,n,3n+i-k,0,l),
\end{align*}
where by setting $u=4n-a$ and $v=3n+a+b-c-j-l$
\begin{align*}
D(a,b,c,j,l)=&\sum_{s=0}^{u} {u+v-s-1\choose v-1} \frac{E(c,l,s)}{(-1)^{v}}
+\sum_{s=0}^{v} {u+v-s-1\choose u-1} \frac{E'(c,l,s)}{(-1)^{v-s}},\\
D'(i,k,l)=&
\sum_{s=0}^{4n+k-i-l} {7n-l-s-1\choose 3n+i-k-1} \frac{E''(l,s)}{ 2^{7n-l-s}}
+\sum_{s=0}^{3n+i-k} {7n-l-s-1\choose 4n+k-i-l-1} \frac{E''(l,s)}{ 2^{7n-l-s}}.
\end{align*}
Here by setting $\gb_{s,t}=-B_sB_tB_{9n-s-t}/(s!t!(9n-s-t)!)$,
$\gb''_{a,b}=\gb_{a,b}(1+(2^{1-a}-1)(2^{1-b}-1))/2$, and
$\gb'''_{a,b,c}=\gb_{a,b}(1+(2^{1-a}-1)(2^{1-b}-1)(2^{1-c}-1))/2$,  we have
{\allowdisplaybreaks
\begin{align*}
E(c,l,s)=&\sum_{t=0}^c {c+s-t-1\choose s-1} \frac{\gb_{l,t}}{2^{c+s-t}}
+\sum_{t=0}^s {c+s-t-1\choose c-1}\frac{\gb_{l,t}}{2^{c+s-t}},\\
E'(c,l,s)=&\sum_{t=0}^c {c+s-t-1\choose s-1} \gb_{l,t}
+\sum_{t=0}^s {c+s-t-1\choose c-1} \gb_{l,t},\\
E''(l,s)=&
\sum_{t=0}^n {n+s-t-1\choose s-1} \gb'''_{l,t,9n-l-t}
+\sum_{t=0}^s {n+s-t-1\choose n-1}  \gb''_{l,t}.
\end{align*}}
\end{thm}
For example, Maple computation shows that
\begin{align*}
\gz_W(2,\soseven)
=&\frac{2^3\cdot 19}{3^3\cdot 7\cdot 17!}\pi^{18},\\
\gz_W(4,\soseven)
=&\frac{2^{12}\cdot 307\cdot 267743941589}{3\cdot 7\cdot 13\cdot 19\cdot 37!}\pi^{36},\\
\gz_W(6,\soseven)
=&\frac{2^{21}\cdot 2053\cdot 9079132487\cdot 265178091767}{3\cdot 7\cdot 11\cdot 19\cdot 54!}\pi^{54},\\
\gz_W(8,\soseven)
=&\frac{2^{29}\cdot 241\cdot 40670746903\cdot 36209034431567319455922705846157}{3\cdot 5\cdot 7\cdot 13\cdot 19\cdot 74!}\pi^{72},\\
\gz_W(10,\soseven)
=&\frac{2^{37}\cdot 61\cdot 45197\cdot 3920899\cdot 3246046224154033\cdot 202097025268393295809502658929}{3^3\cdot 7\cdot 11\cdot 19\cdot 31\cdot 89!}\pi^{90}.
\end{align*}
When taking the sum over the range of $|m_1|,|m_2|,|m_3|\le 100$ in \eqref{equ:definso7}
we find that $\gz_W(2,\soseven)$ is correct up to at least 19 digits,
$\gz_W(4,\soseven)$ up to 42 digits, $\gz_W(6,\soseven)$ up to 64 digits,
$\gz_W(8,\soseven)$ and $\gz_W(10,\soseven)$ up to at least 80 digits.

\section{The $\slfive$ case}
By definition
\begin{equation}
\label{equ:definso7}
   \zeta_W(n;\slfive)=288^n   \sum_{m_1,\dots,m_4=1}^\infty
   \left(\frac{ 1/\big(m_1 m_2 m_3 m_4 (m_1+m_2)(m_2+m_3)(m_3+m_4)\big)}
     {(m_1+m_2+m_3) (m_2+m_3+m_4) (m_1+m_2+m_3+m_4) } \right)^{n}.
\end{equation}
The corresponding matrix is
\begin{equation*}
\underset{\ \phantom{\gs=} \ccic{1}\ \hskip2.1ex \ccic{2}\ \hskip2.1ex \ccic{3}\ \hskip2.1ex \ccic{4}\ \hskip2.1ex \ccic{5}\ \hskip2.1ex \ccic{6}\ \hskip2.1ex \ccic{7}\ \hskip2.1ex \ccic{8}\ \hskip2.1ex \ccic{9} \ \hskip2.1ex \ccic{\x} }{ \gs=\left(\begin{matrix}
1_{\phantom n}& 0_{\phantom n}& 0_{\phantom n}& 0_{\phantom n}& 1_{\phantom n}& 0_{\phantom n}& 0_{\phantom n}& 1_{\phantom n}& 0_{\phantom n}& 1_{\phantom n}\\
0_{\phantom n}& 1_{\phantom n}& 0_{\phantom n}& 0_{\phantom n}& 1_{\phantom n}& 1_{\phantom n}& 0_{\phantom n}& 1_{\phantom n}& 1_{\phantom n}& 1_{\phantom n}\\
0_{\phantom n}& 0_{\phantom n}& 1_{\phantom n}& 0_{\phantom n}& 0_{\phantom n}& 1_{\phantom n}& 1_{\phantom n}& 1_{\phantom n}& 1_{\phantom n}& 1_{\phantom n}\\
0_n& 0_n& 0_n& 1_n& 0_n& 0_n& 1_n& 0_n& 1_n& 1_n
\end{matrix}\right).}
\end{equation*}
Here ${}^{ \bigcirc\hskip-1.25ex{\x}\hskip.3ex{}}$ is the 10th column. The set of 3-column dependencies is
$$\{(1, 2, 5), (1, 6, 8), (1, 9, 10), (2, 3, 6), (2, 7, 9), (3, 4, 7), (3, 5, 8), (4, 6, 9), (4, 8, 10), (5, 7, 10)\}.$$
So we have the following computation tree
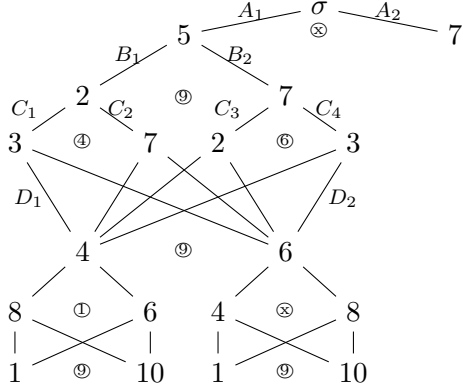
\begin{center}
\begin{tikzpicture}[scale=0.9]
\node (gs) at (-1.5,0.8) {$\gs$};
\node (A1) at (-3.5,0.4) {5};
\node (A2) at (0.5,0.4)  {7};
\draw (-1.5,0.45)   node  {$\ccic{\x}$};
\node (B1) at (-5,-0.5) {2};
\node (B2) at (-2,-0.5) {7};
\draw (-3.5,-0.55)   node {$\ccic{9}$};
\node (C1) at (-6,-1.2) {3};
\node (C2) at (-4,-1.2) {7};
\draw (-5,-1.2)  node  {$\ccic{4}$};
\node (C3) at (-3,-1.2) {2};
\node (C4) at (-1,-1.2) {3};
\draw (-2,-1.2)  node {$\ccic{6}$};

\node (D1) at (-5,-2.8) {4};
\node (D2) at (-2,-2.8) {6};
\draw (-3.5,-2.8) node {$\ccic{9}$};

\node (E1) at (-6,-3.7) {8};
\node (E2) at (-4,-3.7) {6};
\draw (-5,-3.7) node {$\ccic{1}$};
\node (E3) at (-3,-3.7) {4};
\node (E4) at (-1,-3.7) {8};
\draw (-2,-3.7) node {$\ccic{\x}$};

\node (F1) at (-6,-4.6) {1};
\node (F2) at (-4,-4.6) {10};
\draw (-5,-4.6) node {$\ccic{9}$};
\node (F3) at (-3,-4.6) {1};
\node (F4) at (-1,-4.6) {10};
\draw (-2,-4.6) node {$\ccic{9}$};

\draw (gs) to (A1) node[midway] {${\scriptstyle A_1}\atop\ $};
\draw (gs) to (A2) node[midway] {${\scriptstyle A_2}\atop\ $};
\draw (A1) to (B1) node[midway] {${{\scriptstyle B_1}\atop }\ $};
\draw (A1) to (B2) node[midway] {$\ {\scriptstyle B_2}\atop $};
\draw (B1) to (C1) node[midway] {${{\scriptstyle C_1}\atop }\quad\ \ $};
\draw (B1) to (C2) node[midway] {$\ {\scriptstyle C_2}\atop $};
\draw (B2) to (C3) node[midway] {${\scriptstyle C_3 \atop \ }\quad\ \ $};
\draw (B2) to (C4) node[midway] {$\ \ {\scriptstyle C_4\atop \ }$};
\draw (C1) to (D1) node[midway] {${}_{\ \atop\scriptstyle D_1  }\quad\ $};
\draw (C1) to (D2) ;
\draw (C2) to (D1) ;
\draw (C2) to (D2) ;
\draw (C3) to (D1) ;
\draw (C3) to (D2) ;
\draw (C4) to (D1) ;
\draw (C4) to (D2) node[midway] {$\quad\ \ {}_{\  \atop \scriptstyle D_2} $};
\draw (E1) to (D1);
\draw (E2) to (D1);
\draw (E3) to (D2);
\draw (E4) to (D2);
\draw (E1) to (F1);
\draw (E1) to (F2);
\draw (E2) to (F1);
\draw (E2) to (F2);
\draw (E3) to (F3);
\draw (E3) to (F4);
\draw (E4) to (F3);
\draw (E4) to (F4);
\end{tikzpicture}
\end{center}
By symmetry $A_1=A_2$ so we get
\begin{thm}
Let $n$ be a positive even integer. Set ${t \choose -1}=0$ for all $t$. Define
$\gb_{s,t,k}=0$ if $s=1$ or $t=1$ or $k=1$ and define
$\gb_{s,t,k}=(2\pi)^{10n} B_sB_tB_kB_{w-s-t-k}/(s!t!k!(w-s-t-k)!)$
for all other nonnegative integers $s,t,k$. Then
$$\frac{60\gz_W(n,\slfive)}{288^n}=\sum_{i=0}^n{2n-i-1\choose n-1} \left(
 \sum_{j=0}^i {n+i-j-1\choose n-1} B_1(i,j)
+\sum_{j=0}^n {n+i-j-1\choose i-1} B_2(i,j)\right)$$
where for $\ga=1,2$
\begin{align*}
B_\ga(i,j)=&
 \sum_{k=0}^j {n+j-k-1\choose n-1} C_\ga(i,j,k)
+\sum_{k=0}^n {n+j-k-1\choose j-1} C_\ga(i,j,k),\\
C_\ga(i,j,k)=&
\sum_{l=0}^n {3n+j-k-l-1\choose 2n+j-k-1} D_\ga(i,k,l)
+ \sum_{l=0}^{2n+j-k} {3n+j-k-l-1\choose n-1} D_{3-\ga} (i,k,l),\\
D_\ga(i,k,l)=&\sum_{s=0}^l {n+l-s-1\choose n-1} E_\ga(i,k,l,s)
+ \sum_{s=0}^n {n+l-s-1\choose l-1} E_\ga(i,k,l,s),
\end{align*}
and
\begin{align*}
E_1(i,k,l,s)=& \sum_{t=0}^{2n+l-s} {5n+l-i-s-t-1\choose 3n-i-1}\gb_{s,t,k}
+\sum_{t=0}^{3n-i} {5n+l-i-s-t-1\choose 2n+l-s-1}\gb_{s,t,k},\\
E_2(i,k,l,s)=& \sum_{t=0}^n {5n+l-i-s-t-1\choose 4n+l-i-s-1} \gb_{s,t,k}
+ \sum_{t=0}^{4n+l-i-s} {5n+l-i-s-t-1\choose n-1} \gb_{s,t,k}.
\end{align*}
\end{thm}
For example, we have
\begin{align*}
\gz_W(2,\slfive)
=&\frac{1}{650970015609375}\pi^{20}=\frac{2^{16}\cdot 13}{3^2\cdot 5^3\cdot 7\cdot 11\cdot 18!}\pi^{20},\\
\gz_W(4,\slfive)
=&\frac{2^{38}\cdot 1523\cdot 2625375581}{3^2\cdot 5^2\cdot 7\cdot 11\cdot 41!} \pi^{40},\\
\gz_W(6,\slfive)
=&\frac{2^{57}\cdot  30677\cdot  2082905565627654787323001}{3^2\cdot 5^2\cdot 7\cdot 11\cdot 13\cdot 31\cdot  61!}\pi^{60},\\
\gz_W(8,\slfive)
=&\frac{2^{79}\cdot 3^2\cdot 11\cdot 85081\cdot 1361779882876127669651\cdot 728520415874861 }{ 5^2\cdot 7\cdot 17\cdot 82!}\pi^{80},\\
\gz_W(10,\slfive)
=&  \frac{2^{98}\cdot 29\cdot 13^2\cdot 2143\cdot 6751027\cdot 430667831149 a }{3^2\cdot 5^2\cdot 11\cdot 17\cdot 101!}  \pi^{100},
\end{align*}
where $a=201223346979560452521803194127591413.$

\end{document}